\documentclass[12pt,reqno]{amsart}
\usepackage{amsmath,amssymb, amsthm}

\usepackage[colorlinks,urlcolor=blue]{hyperref}
\usepackage{graphicx}
\usepackage{multirow}
\usepackage{cite}
\usepackage{algpseudocode}
\usepackage{algorithm}
\usepackage{array}
\newcolumntype{L}[1]{>{\raggedright\let\newline\\\arraybackslash\hspace{0pt}}m{#1}}
\newcolumntype{C}[1]{>{\centering\let\newline\\\arraybackslash\hspace{0pt}}m{#1}}
\newcolumntype{R}[1]{>{\raggedleft\let\newline\\\arraybackslash\hspace{0pt}}m{#1}}

\theoremstyle{plain}
\newtheorem{thm}{Theorem}[section]

\newtheorem{prop}[thm]{Proposition}
\numberwithin{equation}{section} \theoremstyle{definition}

 \newtheorem{exmp}[thm]{Example}

\allowdisplaybreaks

\numberwithin{equation}{section}

\renewcommand{\to}{\longrightarrow}
\begin{document}
\title[A practical approach to optimization]
{A practical approach to optimization}
\author[W. Jirakitpuwapat et al.]
{Sompong Dhompongsa$^{1,2}$, Wachirapong Jirakitpuwapat$^{1}$,  Konrawut Khammahawong$^{2}$ and Poom Kumam$^{1, 2, \S}$, 
}
\date{}
\thanks{\it $^\S$ Corresponding author: poom.kum@kmutt.ac.th (P. Kumam)}
\maketitle

\begin{flushleft}
	{\footnotesize $^{1}$KMUTTFixed Point Research Laboratory, KMUTT-Fixed Point Theory and Applications Research Group, SCL 802 Fixed Point Laboratory, Department of Mathematics, Faculty of Science, King Mongkut's University of Technology Thonburi (KMUTT), 126 Pracha-Uthit Road, Bang Mod, Thrung Khru, Bangkok 10140, Thailand
		\vskip 2mm
		$^{2}$Center of Excellence in Theoretical and Computational Science (TaCS-CoE), Science Laboratory Building,
		King Mongkut's University of Technology Thonburi (KMUTT), 126 Pracha-Uthit Road, Bang Mod, Thrung Khru, Bangkok 10140, Thailand 
		\vskip 2mm
		Email addresses: sompong.dho@kmutt.ac.th (S. Dhompongsa) wachirapong.jira@hotmail.com (W. Jirakitpuwapat),  and k.konrawut@gmail.com (K. Khammahawong) and poom.kum@kmutt.ac.th  (P. Kumam), }
\end{flushleft}

\maketitle

\hrulefill

\begin{abstract}
We present a new approach for finding a minimal value of an arbitrary function assuming only its continuity. The process avoids verifying Lagrange- or KKT-conditions. The method enables us to obtain a Brouwer fixed point (of a continuous function mapping from a cube into itself).
\end{abstract}

\begin{footnotesize}



\end{footnotesize} 	
\maketitle

\hrulefill

\section{Introduction}
For a given set of continuous functions 
$f,g_1,g_2,\dots,g_m,h_1,h_2,\dots,h_n:\mathcal{C}=\prod_{i=1}^{p}
[c_i,d_i] $ $\to \mathbb{R}$,
a minimization problem of the form
\begin{equation} \label{main}
\begin{aligned}
& \min\limits_{x \in \mathcal{C}}f(x) \\
\text{subject to } & g_i(x) = 0 \text{ } (i=1,2,\dots,m)\\
& h_i(x) \leq 0 \text{ } (i=1,2,\dots,n).\\
\end{aligned}
\end{equation}
is well known. For the Problem \eqref{main}, $f$ is called the objective function and the equalities (described by $g_i$) and the inequalities (described by $h_i$) are called the constraints. We call the set $\mathcal{A} = \{x \in \mathcal{C}:g_i(x)=0 \text{ } (i=1,2,\dots,m) \text{ and } h_i(x)\leq 0\text{ } (i=1,2,\dots,n)\}$ the feasible set of Problem \eqref{main}. If $\mathcal{A}$ is not empty, it is compact since it is a zero set of the continuous function $F$ defined below.  Consequently, Problem \eqref{main} always has a solution if $\mathcal{A}$ is not empty.   The subject is well understood for convex optimization with Lagrange multipliers and
Karush-Kuhn-Tucker conditions are its familiar main tools. It is the purpose of this article to introduce an alternative method in minimizing a function without using the tools mentioned above.  The method transforms the constrained Problem \eqref{main} of $f$ into an unconstrained one of a deformation $f_t$ of $f.$ It can be considered as a toolkit using for approximating a result by applying any existing software. We choose to work on some well-known software to find a decreasing sequence $\{f_t (x_n)\}$, namely, particle swarm optimization (PSO), particle-search algorithm, and convex optimization. By testing the method over many kinds of objective functions $f,$  we believe the method is quite practical. It is found that a problem may work well under one software but not under some others. Moreover, the method can be performed to obtain a Brouwer fixed point and  applied to a vector optimization.

In computational science, particle swarm optimization (PSO) \cite{Kennedy,Mezura_Montes_2011,pedersen2010good} is the  computational method that optimization  problem by iteratively trying to improve a candidate solution with regard to a given measure of quality.  A basic variant of the PSO algorithm works by having a population (swarm) of candidate solutions (particles). These particles are moved around in the search-space according to a simple formula. The movements of the particles are guided by their own best known position in the search-space. The entire swarm's best known position.  When improved positions are being discovered these will then come to guide the movements of the swarm. The process is repeated and by doing so it is hoped, but not guaranteed, that a satisfactory solution will eventually be discovered.

Pattern search algorithm is a family of numerical optimization methods. It finds a sequence of points that approach an optimal point. The value of the objective function either decreases or remains the same from each point in the sequence to the next  \cite{abramson2002pattern,MR1972220,Conn_1997}.

Convex optimization is a subfield of mathematical optimization that studies the problem of minimizing convex functions over convex sets.  Convex algorithm is a mathematical method of solving convex optimization  \cite{MR1795061,MR1724768,MR1387333}. The key to the algorithmic success in minimizing convex functions is that these functions exhibit a local to global phenomenon. This local to global phenomenon is that local minimal of convex functions are in fact global minimal.

\section{Methodology}

Put $G_i=|g_i| \text{ }$ $(i=1,2,\dots,m)$, $H_i=|h_i|+h_i \text{ }$ $(i=1,2,\dots,n)$, and $F = \sum\limits_{i=1}^{m}G_i+\sum\limits_{i=1}^{n}H_i$. 
Clearly, $F$ is continuous and $F(x)=0$ if and only if $x$ satisfies the constraints of Problem \eqref{main} (i.e., it lies in the feasible set $\mathcal{A}$). For large numbers $K$ and $M$, set for $t \in (0,1), f_t=(1-t)(f-K)+tMF $.\\
Since we are going to work on the deformed function $f_t$ for $t$ sufficiently close to $1$, we therefore take any existing software available. We select $3$ softwares, namely Particle Swarm Optimization, Pattern-Search, and Convex Algorithm. We let $K$ to be large to be certained that the graph of $f-K$ totally lies under the graph of $F.$     As for large $M$, we try to make it easy for a software to find a decreasing sequence $\{f_t (x_n)\}$.  The parameter $t$ getting close to $1$ is to making the iteration point $x_n$ being closer to or lying in the feasible set $\mathcal{A}.$
\begin{prop}\label{prop 3.1}
	For any $t \in (0,1)$ with $f_t >0$ outside $\mathcal{A},$ $x$ is a minimizer of Problem \eqref{main} if and only if $x$ is a minimizer of $f_t$.
\end{prop}
\begin{proof}
	This is straightforward since $f_t=(1-t)(f-K)$ on $\mathcal{A}$.
\end{proof}
\noindent By the term ``minimizer'' it is meant to be a minimal element, i.e., a local minimizer.
\begin{algorithm}[h] 
	\caption{Example code (PAO our Algorithm)}
	\begin{algorithmic}
		\State \textbf{Input} Set up problem \ref{main}
		\State \textbf{Parameter} $K,M,t$
		\State \textbf{Output} $x$ \\
		\qquad	$G_i=|g_i| \ $ $(i=1,2,\dots,m)$\\
		\qquad	$H_i=|h_i|+h_i \ $ $(i=1,2,\dots,n)$\\
		\qquad	$F = \sum\limits_{i=1}^{m}G_i+\sum\limits_{i=1}^{n}H_i$\\
		\qquad	$f_t=(1-t)(f-K)+tMF$\\
		\qquad	$x = \arg \min\limits_{x \in C}f_t(x)$ \\
	\end{algorithmic} 	
	\label{alg}
\end{algorithm}
\section{Applications}

\subsection{Brouwer Fixed Points}
The Brouwer fixed theorem says that any continuous mapping $T=(f_1,\dots,f_d):\prod_{i=1}^{d}[a_i,b_i] \rightarrow \prod_{i=1}^{d}[a_i,b_i]$ always has a fixed point. See \cite{MR1511679,MR3745100,MR4009295,MR3446192} for some new proofs. To find a fixed point of $T$, set in Problem \eqref{main}, $f(x_1,x_2,\dots,x_d)=1$ and $g_i(x_1,x_2,\dots,x_d)=f_i(x_1,x_2,\dots,x_d)-x_i$ $(i=1,2,\dots,d)$. (See Example \ref{ex 5.5} and \ref{ex 5.6}.)

\subsection{Vector Optimization}
Given continuous mappings $f_1,f_2,\dots,f_k,g_1,g_2,\dots,g_m,h_1,$ $h_2,\dots,h_n:\mathcal{C}=\prod_{i=1}^{p}
[c_i,d_i] $ $\to \mathbb{R}$.
We need to solve
\begin{equation} \label{VectorOptimization}
\begin{aligned}
& \min\limits_{x \in \mathcal{C}}(f_1(x),f_2(x),\dots,f_k(x)) \text{ (with respect to an order)}  \\
\text{subject to } & g_i(x) = 0 \text{ } (i=1,2,\dots,m)\\
& h_i(x) \leq 0 \text{ } (i=1,2,\dots,n).\\
\end{aligned}
\end{equation}
We consider the problem of the forms: 
\begin{itemize}
	\item[(1)] $\min\limits_{x \in \mathcal{C}}\sum\limits_{i=1}^{k}f_i(x)$. Set  $f = \sum_{i=1}^{k}f_i$ for the objective function in Problem \eqref{main}. (See Example \ref{ex 5.8}.)
	\item[(2)] Finding $x^\ast=(x_1^\ast,x_2^\ast,\dots,x_p^\ast) \in \mathcal{C}$ such that $f_i(x^\ast) \leq c_i$, where $c_i\leq t_i$ for some thresholds $t_i$ \text{ } $(i=1,2,\dots,k)$. To comply with Problem \eqref{main}, we set $f=1$ as an objective function and additionally define $h_i=f_i-c_i$ \text{ } $(i=n+1,n+2,\dots,n+k)$. (See Example \ref{ex 5.9}.)
\end{itemize}
In practice, if we only want to find a point $x^*$ with $f(x^*) \leq c$ for some assigned number $c$, Problem \eqref{main} can read as 
\begin{equation} 
\begin{aligned}
& \min\limits_{x \in \mathcal{C}}1 \\
\text{subject to } & g_i(x) = 0 \text{ } (i=1,2,\dots,m)\\
& h_i(x) \leq 0 \text{ } (i=1,2,\dots,n)\\
&f(x) - c\leq 0. \\
\end{aligned}
\end{equation}

\section{Numerical Examples}
We choose $\mathcal{C} = [-10,10]^p, K = 100, M = 10000$ and $t=0.95$. We experiment on nine Examples, and record results in three Tables. 
The Tables display approximate  minimizers and constraint validation.

\begin{exmp}\label{ex 5.1} \cite{Franklin_1980}
$$
	\begin{array}{ll} 
	\min\limits_{x \in \mathcal{C}} & x_1^2+x_1x_2+x_2^2-5x_2 \\
	\text{subject to } & x_1+x_2 = 1 \\ 
	& x_1 \geq 0 \\ 
	& x_2 \geq 0 \\
	\end{array}
$$
\end{exmp}

\begin{exmp}\label{ex 5.2} \cite{Franklin_1980}
$$
	\begin{array}{ll}
	\min\limits_{x \in \mathcal{C}} & -(x_1-3)^6-(x_2-4)^6 \\
	\text{subject to } & x_1^2+x_2^2 \leq 25 \\
	& x_1+x_2 \geq 7 \\
	& x_1 \geq 0 \\
	& x_2 \geq 0 \\
	\end{array}
$$
\end{exmp}

\begin{exmp}\label{ex 5.3}\cite{Franklin_1980}[Geometric Programming]
$$
		\begin{array}{lll}
	\min\limits_{x \in \mathcal{C}} & \frac{1}{x_1x_2x_3}+x_1x_2\\
	\text{subject to } & 0.5x_1x_3+0.25x_1x_2 \leq 1\\
	& x_1 \geq 0 \\
	& x_2 \geq 0 \\
	& x_3 \geq 0 \\
	\end{array}
$$
\end{exmp}

\begin{exmp}\label{ex 5.4}\cite{Franklin_1980}
$$
	\begin{array}{lll}
	\min\limits_{x \in \mathcal{C}} & \frac{1}{x_1x_2x_3}+x_1x_2+x_3^7\\
	\text{subject to } & 0.5x_1x_3+0.25x_1x_2 \leq 1\\
	& x_1 \geq 0 \\
	& x_2 \geq 0 \\
	& x_3 \geq 0 \\
	\end{array}
$$
\end{exmp}

\begin{exmp}\label{ex 5.45}
$$
	\begin{array}{ll}
	\min\limits_{x \in \mathcal{C}} & 4x_1+10x_2+15x_3\\
	\text{subject to } & x_1+2x_2+3x_3 = 3\\
	& 3x_1+x_2+2x_3 = 7.5\\
	& x_1 \geq 0 \\
	& x_2 \geq 0 \\
	& x_3 \geq 0 \\
	\end{array}
$$
\end{exmp}

\begin{exmp}\label{ex 5.5}
$$
	\begin{array}{ll}
	\min\limits_{x \in \mathcal{C}} & 1\\
	\text{subject to } & 0.5(\cos(x_1+x_2-x_3^4x_5))x_4-x_1 = 0 \\
	& 0.1(|x_1x_2+x_3-x_5|+x_4^2)-x_2 = 0 \\
	& (x_1+x_3x_4-(x_2+x_5)^2)/30-x_3 = 0\\
	& (x_1-x_2^2+x_3-x_5^2)/12-x_4  = 0\\
	& (x_1+x_2-(x_3+x_5+x_4)^2)/40-x_5 = 0\\
	\end{array}
$$
\end{exmp}

\begin{exmp}\label{ex 5.6}
$$
	\begin{array}{lll}
	\min\limits_{x \in \mathcal{C}} & 1\\
	\text{subject to } & 0.001((x_1+3)^2+(x_2-2)^4+x_3^2+x_4^2+x_5)-x_1 = 0 \\
	& 0.01(x_1+(x_2+5)^2+x_3+x_4+(x_5+2))-x_2 = 0 \\
	& 0.001(x_1^4+(x_4-3)^2+(x_5+2)^2)-x_3 = 0\\
	& 0.001((x_3-3)^4+x_5^2+x_1^4)-1-x_4  = 0\\
	& 0.01(x_1^2+x_2+x_3-(x_5-1)^2)-x_5 = 0\\
	\end{array}
$$
\end{exmp}

\begin{exmp}\label{ex 5.8}\cite{Franklin_1980}
$$
	\begin{array}{lll}
	\min\limits_{x\in \mathcal{C}} & (x_1^2 - 5x_1 + 7x_2) + (-x_1^2 -x_2^2) + (x_1-1)^2 +(x_2-5)^2 \\
	\text{subject to } & 3x_1 + 4x_2 =6\\
	& x_1+x_2 =2\\
	& 2x_1 + 3x_2 \leq 6 \\
	& x_1 \geq 0 \\
	&x_2 \geq 0
	\end{array}
$$
\end{exmp}
\begin{exmp}\label{ex 5.9}\cite{Franklin_1980}
$$
	\begin{array}{lll}
	&4x_1^2+x_2^2-x_1-2 \leq 1\\
	&e^{-x_1}-x_1-2x_2 \leq 1 \\
	\text{subject to } &  2x_1+x_2 \leq 1\\
	& x_1^2\leq 1 \\
	& \sqrt{x_1^2+x_2^2 }-x_1^3 \leq 2 \\
	& -x_1^3+0.5(-x_2-x_2^3+|x_2^3-x_2|) \leq 0, \quad x_1,x_2 \in \mathbb{R}.                \\
	\end{array}
	$$
\end{exmp}

{\footnotesize
	\begin{table}[h] 
	\centering
	\caption{Particle Swarm Optimization}\label{TablePSO}
	\begin{tabular}{|c|c|c|C{5.5cm}|c|c|}
		\hline
		\multirow{2}{*}{Example}           & \multicolumn{5}{c|}{PSO}                                                                                                                                                                                                                                                                                                                                                                                                                                                                                                                                               \\ \cline{2-6}
		&initial point            & value                                                    & $x$                                                                                                                                                                                                                                                                                            & $\max\limits_{x \in \mathcal{C}}|g_i(x)|$ & $\max\limits_{x \in \mathcal{C}}h_j(x)$  \\ \hline
		$\ref{ex 5.1}$ &-   & $-4$ & $(0,1)$                                                                                                                                                                                                                                                                           & $0$                                                                                                         & $0$                                                                                                   \\ 
		$\ref{ex 5.2}$  &-  & $-2$                                                   & $(4,3)$                                                                                                                                                                                                                                                                              & $-$                                                                                                         & $0$                                                                                      \\
		$\ref{ex 5.3}$   &-  & $0.6325$                                                   & 
		$(10,0.0316,10)$                                                                                                                                                                                                                                                                     & -                                                                                      & $-0.0316 $                                                                                                  \\ 
		$\ref{ex 5.4}$  &- & $2.4397$                                                   & $(0.1141,9.9975,0.7715)$                                                                                                                                                                                                                          & - & $-0.1141$                                                                                                         \\  
		$\ref{ex 5.45}$  &- & $12.6$                                                   & $(2.4,0.3,0)$                                                                                                                                                                                                                          & $0$ & $0$                                                                                                         \\ 
		$\ref{ex 5.5}$ &-  & $1$                                                        & 
		$(-1.977 \times 10^{-11},1.02 \times 10^{-12},-1.067 \times 10^{-12},-2.719 \times 10^{-11},  6.04 \times 10^{-12})$ & $2.546 \times 10^{-11} $                                                     & $-$                                                                                                         \\ 
		$\ref{ex 5.6}$  &-   & $1$                                                        & 
		$(0.018,0.291,0.019,-0.921,-0.007)$                                                                                                       & $1.766\times 10^{-12}$                                                     & $-$                                                                                                        \\                     
		$\ref{ex 5.8}$  &-   & $16$                                                        & 
		$(2,0)$                                                                                                       & $0$                                                     & $0$                                                                                                        \\      
		$\ref{ex 5.9}$  &-   & $1$                                                        & 
		$ (0.7312,1.0271)$                                                                                                       & $-$                                                     & $ -0.4654$                                                                                                        \\       \hline                                                                                           
		
	\end{tabular}
\end{table}
}

{\footnotesize
\begin{table}[h] \label{TablePattern-Search}
	\caption{Pattern-Search Optimization}
	\centering
	\begin{tabular}{|c|c|c|c|c|c|}
		\hline
		\multirow{2}{*}{Example}            & \multicolumn{5}{c|}{Pattern-Search}                                                                                                                                                                                                                                                                                                                                                                                        \\ \cline{2-6}
		&initial point                 & value                                                    & $x$                                                                                                                                                                                                                                                                                            & $\max\limits_{x \in \mathcal{C}}|g_i(x)|$ & $\max\limits_{x \in \mathcal{C}}h_j(x)$ 
		\\ \hline
		$\ref{ex 5.1}$  &$(1,1)$     & $-4$       & $(0,1)$                                                                                                                                                                                  & $0$                                                                                                         & $0$                                                                                                        \\ 
		$\ref{ex 5.2}$ &$(1,1)$    & $-94.3669$ & $(4.6094,1.9374)$                                                                                                                                                                             & $-$                                                                                                         & $0.4532$                                                                                                        \\ 
		$\ref{ex 5.3}$  &$(1,1,1)$   & $0.6325$       & $(0.6325,0.5,10)$                                                                                                                                                                                & $-$                                                                                                         & $-0.5$                                                                                                        \\ 
		$\ref{ex 5.4}$ &$(1,1,1)$  & $2.4397$  & $(1.1385,1,0.7715,1)$                                                                                                                                                                           & $-$                                                                                                    & $-0.2762$                                                                                                         \\  
		$\ref{ex 5.45}$ &$(1,1,1)$  & $12.6429$  & $(2.3571,0,0.2143)$                                                                                                                                                                           & $1.5259  \times 10^{-5}$                                                                                                    & $0$                                                                                                         \\  
		$\ref{ex 5.5}$ &$(0,0,0,0,0)$  & $1$       & $(0,0,0,0,0)$                                                                                                                                                                            & $0$                                                                                                         & $-$                                                                                                         \\
		$\ref{ex 5.6}$ &$(0,0,0,0,0)$    & $1$       & $(0.019,0.029,1.93,-0.92,-0.007)$ & $3.978\times 10^{-6}$                                                     & $-$                                                                                                         \\            
		$\ref{ex 5.8}$  &(1,1)   & $18.7777$                                                        & 
		$(0.6667,1)$                                                                                                       & $1.5259 \time 10^{-5}$                                                     & $-0.6667$                   \\                            
		$\ref{ex 5.9}$  &(1,1)   & $1$                                                        & 
		$(0,1)$                                                                                                       & $-$                                                     & $-1$                   \\ \hline                                                                 
	\end{tabular}
\end{table}
}

{\footnotesize
\begin{table}[h] \label{TableConvexAlgorithm}
	\centering
	\caption{Convex Algorithm}
	\begin{tabular}{|c|c|c|C{5.5cm}|c|c|}
		\hline
		\multirow{2}{*}{Example}   & \multicolumn{5}{c|}{Convex Algorithm}                                                                                                                                                                                                                                                                                                                                                                                                                                                                                                                                                 \\ \cline{2-6}
		&initial point                     & value                                                    & $x$                                                                                                                                                                                                                                                                                            & $\max\limits_{x \in \mathcal{C}}|g_i(x)|$ & $\max\limits_{x \in \mathcal{C}}h_j(x)$  
		\\ 
		$\ref{ex 5.1}$  &$(1,1)$     & $-3.9694$ & $(0.0076,0.9924)$                                                                                                                                                                                                                                                                           & $7.3 \times 10^{-9}$                                                                                                         & $-0.0076$                                                                                                  \\ 
		$\ref{ex 5.2}$  &$(1,1)$   & $-1.2957$                                                   & $(3.9302,3.0698)$                                                                                                                                                                                                                                                                           & $-$                                                                                                         & $-7.97 \times 10 ^{-13}$                                                \\
		$\ref{ex 5.3}$   &$(1,1,1)$    & $0.6325$                                                        & $(0.5623,0.5623,10)$                                                                                                                                                                                                                                                                                   & $-$                                                                                                         & $-0.5623$                                                                                                       \\ 
		$\ref{ex 5.4}$  &$(1,1,1)$   & $2.4397$                                         & $(1.0670,1.0670,0.7715)$                                                                                                                                                                                                                                                                    & $-$                                                                                                   & $-0.3038$                                                                                                        \\ 
		$\ref{ex 5.45}$  & $(1,1,1)$   & $12.6392$                                         & $(2.3608,0.0253,0.1962)$                                                                                                                                                                                                                                                                    & $0.3167 \times 10^{-7}$                                                                                                   & $-0.0253$                                                                                                        \\  
		$\ref{ex 5.5}$ &$(0,0,0,0,0)$    & $1$                                                        & 
		$(-1.520 \times 10^{-10},1.283 \times 10^{-11},1.265 \times 10^{-10},2.282 \times 10^{-10},-3.341 \times 10^{-11})$ & $2.661 \times 10^{-10}$                                                                         & $-$                                                                                                        \\ 
		$\ref{ex 5.6}$	   &$(0,0,0,0,0)$ & $1$                                                        & 
		$(0.018,0.291,0.019,-0.921,-0.007)$                                                                                                   & $8.413 \times 10^{-9}$                                                     & $-$ \\          
		$\ref{ex 5.8}$  &(1,1)   & $16$                                                        & 
		$(2,0)$                                                                                                       & $0$                                                     & $0$                   \\    
		$\ref{ex 5.9}$  &(1,1)   & $1$                                                        & 
		$ (-0.0888,0.8020) $                                                                                                       & $-$                                                     & $ -0.8013$                   \\ \hline   
	\end{tabular}
\end{table}}

\section{Discussion}
In this paper, we transform a constrained optimization to an unconstrained one. Under our approach, the given objective function $f$ (subjected to some constraints) is replaced by a deformed function $f_t$ (without constraints) for some $t$. We chose to use some software packages to approximate a minimizer of $f_t$. We observe that all outcomes approximately satisfy corresponding constraints. Of course, we may obtain different minimizers from different software. It is challenging to construct a new algorithm for finding a global minimizer even for some special cases.

\section*{Acknowledgements}
The authors acknowledge the financial support provided by the Center of Excellence in Theoretical and Computational Science (TaCS-CoE),
KMUTT. The first author  was supported by the ``Petchra Pra Jom Klao Ph.D. Research Scholarship from King Mongkut's University of Technology Thonburi''  (No. Grand 10/2560). 

\bibliographystyle{vancouver}

\end{document}